\documentclass[10pt,reqno]{amsart}
\usepackage[utf8]{inputenc} % UTF-8 encoding
\usepackage[T1]{fontenc}    % Proper font encoding
\usepackage{lmodern}  
\usepackage{bm}

\usepackage{graphics,color}
\usepackage{amssymb, amsmath, amsthm, amscd}
\usepackage{latexsym, verbatim, graphicx, amsfonts}
\usepackage{hyperref}
\usepackage[mathscr]{euscript}
\usepackage{amsmath, amsthm, amssymb}
\usepackage{dsfont}
\usepackage{enumitem}
\usepackage{mathrsfs}
\usepackage{mathtools}
\usepackage{float}

\theoremstyle{plain}
\newtheorem{theorem}{Theorem}[section]

\newtheorem{lemma}[theorem]{Lemma}
\newtheorem{proposition}[theorem]{Proposition}
\theoremstyle{definition}

\newtheorem{definition}[theorem]{Definition}

\theoremstyle{remark}
\newtheorem{rem}[theorem]{Remark}

\DeclareMathOperator*{\stone}{ult}

\newcommand{\norm}[1]{\left\lVert#1\right\rVert}

\DeclareMathOperator{\JL}{JL}
\DeclareMathOperator{\fs}{fin}

\newcommand{\finset}[1]{\fs( #1 )}

\newcommand{\vertiii}[1]{{\left\vert\kern-0.25ex\left\vert\kern-0.25ex\left\vert #1 
    \right\vert\kern-0.25ex\right\vert\kern-0.25ex\right\vert}}

\font\sstext=ecss1000
\font\sssub=ecss1000 at 7pt
\font\sssubsub=ecss1000 at 5pt

\newfam\ssfam
\textfont\ssfam=\sstext
\scriptfont\ssfam=\sssub
\scriptscriptfont\ssfam=\sssubsub

\subjclass[2020]{
46E15, % Banach spaces of continuous, differentiable or analytic functions, 
46B42, % Banach lattices 
46B03, %Isomorphic theory (including renorming) of Banach spaces 
54G12. %Scattered spaces
}

\keywords{Complemented subspace, $C(K)$-space, $AM$-space, Banach lattice, Banach space, almost disjoint family}

\begin{document}

\title[The class of Banach lattices is not primary]{The class of Banach lattices is not primary}
\author[A. Acuaviva]{Antonio Acuaviva}
\address{School of Mathematical Sciences,
Fylde College,
Lancaster University,
LA1 4YF,
United Kingdom} \email{ahacua@gmail.com}

\date{\today}

\begin{abstract}
Building on a recent construction of Plebanek and Salguero-Alar\-c\'on, which solved the Complemented Subspace Problem for $C(K)$-spaces, and the subsequent work of De Hevia, Martínez-Cervantes, Salguero-Alarc\'on, and Tradacete solving the Complemented Subspace Problem for Banach lattices,  we show that the class of Banach lattices is not primary. Specifically, we exhibit a compact Hausdorff space $L$ such that $C(L) \simeq X \oplus \tilde{X}$ and neither $X$ nor $\tilde{X}$ is isomorphic to a Banach lattice. In particular, it also follows that the class of $C(K)$-spaces is not primary.
\end{abstract}

\maketitle

\bigskip
\section{Introduction and main results}

In a remarkable construction, Plebanek and Salguero-Alarc\'on \cite{plebanek} constructed compact Hausdorff spaces $K$ and $L$ such that $C(L) \simeq C(K) \oplus \mathrm{PS}_2$ and the Banach space $\mathrm{PS}_2$ fails to be isomorphic to a $C(K)$-space, giving a negative solution to the longstanding Complemented Subspace Problem for $C(K)$-spaces. Shortly after, De Hevia, Martínez-Cervantes, Salguero-Alarc\'on, and Tradacete \cite{de2023counterexample} showed that the space $\mathrm{PS}_2$ is not even isomorphic to a Banach lattice, thus providing a negative solution to the Complemented Subspace Problem for Banach lattices. 

Building on these results, we show that the class of Banach lattices is not primary. More precisely, we obtain the following.

\begin{theorem}\label{th: main-lattices}
    There exist a compact Hausdorff space $L$ and Banach spaces $X$ and $\tilde{X}$ such that 
    \begin{equation*}
        C(L) \simeq X \oplus \tilde{X},
    \end{equation*}
    and neither $X$ nor $\tilde{X}$ is isomorphic to a Banach lattice.
\end{theorem}

In particular, this answers a question of De Hevia and Tradacete \cite[Question 4]{de2025complemented}, which originally motivated our interest in this problem. Observe that this also implies that the class of $C(K)$-spaces is not primary.

Finally, in passing, we note that, with natural modifications, the arguments in \cite[Section 5]{de2023counterexample} show that the same conclusion holds for complex Banach lattices.

\bigskip
\section{Organization and sketch of the proof of Theorem \ref{th: main-lattices}}

Our construction builds on the construction of Plebanek and Salguero-Alarc\'on \cite{plebanek}, together with the subsequent work of De Hevia, Martínez-Cervantes, Salguero-Alarc\'on, and Tradacete \cite{de2023counterexample}.

Section \ref{sec: preliminaries} introduces basic terminology and recalls concepts necessary for the construction. This section also builds the framework in which the construction will be performed. Section \ref{sec: technical-results} reformulates some results of Plebanek and Salguero-Alarc\'on in our setting. The reader familiar with their work can safely skim this section and refer back to it as and when needed. Since the proofs of these modified results are essentially the same as in the work of Plebanek and Salguero-Alarc\'on, we either omit them or provide only outlines of the necessary changes, as appropriate. \\

Section \ref{sec: proof-main-theorem} contains the construction of the counterexample. The idea is as follows. Plebanek and Salguero-Alarc\'on constructed compact Hausdorff spaces $K$ and $L$ such that  
\begin{equation*}
    C(L) = C(K) \oplus \mathrm{PS}_2
\end{equation*}
and $\mathrm{PS}_2$ fails to be a isomorphic to a $C(K)$-space.

Their approach relies on a transfinite construction that rules out the existence of a weak$^*$-compact subset of $B_{\mathrm{PS}_2^*}$ which is both norming and free. 

In our setting, we perform two such constructions simultaneously, producing four compact Hausdorff spaces $K$, $\tilde{K}$, $L$, and $\tilde{L}$ such that  
\begin{equation*}
    C(L) = C(K) \oplus Y, \quad C(\tilde{L}) = C(\tilde{K}) \oplus \tilde{Y}.
\end{equation*}  
This simultaneous approach allows us to intertwine the structure of $K$ with $\tilde{Y}$ and that of $\tilde{K}$ with $Y$ so that neither 
\begin{equation}\label{eq: X-and-X}
    X := C(\tilde{K}) \oplus Y \quad\text{nor}\quad \tilde{X} := C(K) \oplus \tilde{Y}
\end{equation}  
is isomorphic to a Banach lattice, but clearly
\begin{equation*}
    X \oplus \tilde{X} = (C(\tilde{K}) \oplus Y) \oplus (C(K) \oplus \tilde{Y}) = C(L) \oplus C(\tilde{L}) = C(L \sqcup \tilde{L}),
\end{equation*}  
where $L \sqcup \tilde{L}$ denotes the topological sum of $L$ and $\tilde{L}$. 

The proof that the Banach spaces $X$ and $\tilde{X}$ defined by \eqref{eq: X-and-X} are not isomorphic to Banach lattices is now essentially the same as the one in \cite[Section 4]{de2023counterexample}. We include the details in Section \ref{sec: lattice} for the reader's convenience.

We stress that the underlying idea remains the same as in Plebanek and Salguero-Alarc\'on’s work: there are sufficiently many possible choices in their construction (in fact, continuum many), which enables us to select ones that prevent $X$ and $\tilde{X}$ from being isomorphic to Banach lattices.

\bigskip
\section{Preliminaries}\label{sec: preliminaries}

\subsection{Johnson-Lindenstrauss spaces, Boolean algebras, and Alexandrov-Urysohn compacta.} \label{subsec: JohnsonLind}

We denote by $\omega = \{0, 1, 2 \dots, \}$ the set of non-negative integers and let $\finset{\omega}$ be the family of all finite subsets of $\omega$. Given $\mathcal{F} \subseteq \mathcal{P}(\omega)$, we let  $[\mathcal{F}]$ denote the smallest Boolean subalgebra of $\mathcal{P}(\omega)$ containing $\mathcal{F}$.

A family $\mathcal{B}$ of infinite subsets of $\omega$ (or of any other countably infinite set) is called \emph{almost disjoint} if for every distinct pair $A, B \in \mathcal{B}$ the intersection $A \cap B$ is finite.  Given an almost disjoint family $\mathcal{B}$, the \emph{Johnson-Lindenstrauss} space $\JL (\mathcal{B})$ is the closed linear span inside of $\ell_\infty$ of the set of characteristic functions
\begin{equation*}
    \{\mathds{1}_n: n \in \omega\} \cup \{\mathds{1}_B: B \in \mathcal{B}\} \cup \{\mathds{1}_\omega\},
\end{equation*}
where $\mathds{1}_n$ stands for $\mathds{1}_{\{n\}}$. Alternatively, let $\mathfrak{B} = [\mathcal{B \cup \finset{\omega}}]$, then $\JL(\mathcal{B})$ is the closed linear span of simple $\mathfrak{B}$-measurable functions,
that is, functions of the form $f = \sum_{i \leq n} r_i \mathds{1}_{B_i}$, where $n \in \omega$, $r_i \in \mathbb{R}$ and $B_i \in \mathfrak{B}$. It is known that $\JL(\mathcal{B})$ is isometrically isomorphic to a $C(K)$-space. The underlying compact Hausdorff space $\stone(\mathfrak{B})$ is the Stone space of the Boolean algebra $\mathfrak{B}$. Equivalently, this space can be explicitly realized as follows. Define $K_\mathcal{B} = \omega \cup \{p_B: B \in \mathcal{B\}} \cup \{\infty\}$, and specify a topology on $K_\mathcal{B}$ by
\begin{itemize}
    \item points in $\omega$ are isolated,
    \item for every $B \in \mathcal{B}$, a neighbourhoods basis $p_B$ is given by the sets $p_B \cup (B \setminus F)$ where $F \in \finset{\omega}$,
    \item $K_\mathcal{B}$ is the one-point compactification of the locally compact Hausdorff space $\omega \cup \{p_B: B \in \mathcal{B\}}$, endowed with the above topology.
\end{itemize}

It is easy to see that this space is a separable, scattered, compact Hausdorff space with empty third Cantor-Bendixson derivative. Observe that the dual of $\JL(\mathcal{B})$ is isometrically isomorphic to the space $M(\mathfrak{B})$ of real-valued finitely additive measures on $\mathfrak{B}$ of finite variation. 

\subsection{The space $M(\mathfrak{B})$.} \label{subsec: description-measures} The measures in $M(\mathfrak{B})$ have a nice description. Since $\JL(\mathcal{B})$ is isometrically isomorphic to $C(K_\mathcal{B})$, where $K_\mathcal{B}$ is scattered, the Riesz Representation Theorem and Rudin's theorem \cite{rudin1957continuous} imply that
\begin{equation*}
    M(\mathfrak{B}) = \JL(\mathcal{B})^* =  C(K_\mathcal{B})^* = \ell_1(K_\mathcal{B});
\end{equation*}
that is, every measure $\nu \in M(\mathfrak{B})$ can be expressed as
\begin{equation*}
    \nu = \sum_{s \in K_{\mathcal{B}}} c_s \delta_s
\end{equation*}
for some scalars $c_s \in \mathbb{R}$, so that
\begin{equation*}
    \norm{\nu}_1 = \sum_{s \in K_\mathcal{B}} |c_s|.
\end{equation*}

In particular, we can decompose $\nu$ as
\begin{equation*}
    \nu = \mu + \overline{\nu}, \hspace{5pt } \text{ where } \hspace{5pt }\mu = \sum_{s \in \omega} c_s \delta_s \hspace{5pt} \text{ and } \hspace{5pt} \overline{\nu} = \sum_{s \in K_\mathcal{B} \setminus \omega} c_s \delta_s.
\end{equation*}
We will call $\mu$ the $\ell_1(\omega)$-part of $\nu$.

For our purposes, we will usually have an enumeration of the almost disjoint family $\mathcal{B} = \{B_\xi: \xi < \mathfrak{c\}}$, in that case, we can further identify
\begin{equation*}
    M(\mathfrak{B}) = \ell_1(\omega) \oplus \ell_1(\mathfrak{c})
\end{equation*} 
in a canonical way. It is worthwhile to pause briefly to expand on this representation of $M(\mathfrak{B})$, as this will play an important role later. 

Observe that given an element $x \in \ell_1(\omega)$, we can naturally consider it as a measure $\mu \in M(\mathcal{P}(\omega))$ by
\begin{equation*}
    \mu(B) = \sum_{s \in B} x(s), \hspace{5pt} B \subseteq \omega.
\end{equation*}
Furthermore, considering the restriction of $\mu$ to the Boolean subalgebra $\mathfrak{B}$ of $\mathcal{P}(\omega)$, we can consider $\mu \in M(\mathfrak{B})$.

\begin{rem}
    Note that this provides a key insight: the $\ell_1(\omega)$-part has, in a sense, a universal character, 
    as it extends uniquely as a measure to \emph{any} Boolean subalgebra of $\mathcal{P}(\omega)$.
\end{rem}

Now, given a measure $\nu \in M(\mathfrak{B})$, we can decompose $\nu = \mu + \overline{\nu}$, where $\mu \in M(\mathfrak{B})$ arises from an element $x$ of $\ell_1(\omega)$ and $\overline{\nu}$ vanishes on finite sets. Moreover, $|\overline{\nu}|(B) \neq 0$ for at most countably many $B \in \mathcal{B}$, so we may identify $\overline{\nu}$ with an element $y \in \ell_1(\mathfrak{c})$ by setting $y(\alpha) = \overline{\nu}(B_\alpha)$ for $\alpha < \mathfrak{c} $.
Conversely, given $y \in \ell_1(\mathfrak{c})$, we may associate to it a measure $\overline{\nu} \in M(\mathfrak{B})$ in this way. Thus, we may represent any $\nu \in M(\mathfrak{B})$ by a pair $(x,y) \in \ell_1(\omega) \oplus \ell_1(\mathfrak{c})$, where $x$ corresponds to $\mu$ and $y$ to $\overline{\nu}$. Conversely, every $(x,y) \in \ell_1(\omega) \oplus \ell_1(\mathfrak{c})$ determines a measure $\nu = \mu + \overline{\nu} \in M(\mathfrak{B})$. \\

We now record a crucial observation. For any ordinal $\xi < \mathfrak{c}$ consider the subalgebra $\mathfrak{B}(\xi)$ of $\mathfrak{B}$ given by
\begin{equation*}
    \mathfrak{B}(\xi) = [\{B_\alpha: \alpha < \xi \} \cup \finset{\omega}] \subseteq \mathfrak{B},
\end{equation*}
and let $\nu \in M(\mathfrak{B}(\xi))$. Then $\nu$ naturally extends to an element of $M(\mathfrak{B})$. As above, we decompose $\nu = \mu + \overline{\nu}$, where $\mu$ is the $\ell_1(\omega)$-part of $\nu$. Observe that $\mu$ naturally defines a measure in $M(\mathfrak{B})$, while $\overline{\nu}$ can be extended to $M(\mathfrak{B})$ by declaring 
$\overline{\nu}(B_\alpha) = 0$ for every $\xi \leq \alpha < \mathfrak{c} $.

\medskip

Conversely, suppose $\nu \in M(\mathfrak{B})$ and decompose $\nu = \mu + \overline{\nu}$, where $\mu$ is the $\ell_1(\omega)$-part of $\nu$. Since $|\overline{\nu}|(B) \neq 0$ for only countably many $B \in \mathfrak{B}$, there exists $\xi < \mathfrak{c}$ such that $|\overline{\nu}|(B_\alpha) = 0$ 
for all $\alpha \geq \xi$. In particular, we can identify $\nu$ with a measure in $M(\mathfrak{B}(\xi))$, corresponding to its restriction to $\mathfrak{B}(\xi)$. Observe that the natural extension of the latter recovers $\nu$.

\subsection{Complemented subspaces.}

As in the original work of Plebanek and Sal\-gue\-ro-Alar\-c\'on, it will be convenient to work with $\omega \times 2$ instead of $\omega$. Every subset $C$ of the form $C = C_0 \times 2$ for some $C_0 \subseteq \omega$ is called a \emph{cylinder}.  For $n \in \omega$ we denote by $c_n$ the cylinder $\{(n, 0), (n,1)\}$. We say that a set $B \subseteq C$ \emph{splits} the cylinder $C$ if $|B \cap c_n| = 1$ whenever $c_n \subseteq C$. Observe that in this case $C \setminus B$ also splits $C$, and we say that $B$ and $C \setminus B$ is a \emph{splitting of $C$}.

Suppose that $\mathcal{A} = \{A_\xi: \xi < \mathfrak{c\}}$ is an almost disjoint family of cylinders in $\omega \times 2$, together with a splitting $B^i_\xi$, $i = 0,1$ of each cylinder $A_\xi$. Observe that then $\mathcal{B} = \{B^0_\xi, B^1_\xi: \xi < \mathfrak{c}\}$ is also an almost disjoint family. In this situation, we will say that the almost disjoint family $\mathcal{B}$ \emph{splits} $\mathcal{A}$. 

In this setting, we can find a $1$-complemented isometric copy of $\JL(\mathcal{A})$ inside $\JL(\mathcal{B})$, see \cite[Section~3]{plebanek}. With a slight abuse of notation, we will continue to denote this copy by $\JL(\mathcal{A})$. More precisely, this copy is the closed linear subspace of $\ell_\infty(\omega \times 2)$ spanned by
\begin{equation*}
    \{\mathds{1}_{c_n} : n \in \omega \} \cup \{\mathds{1}_A : A \in \mathcal{A} \} \cup \{\mathds{1}_{\omega \times 2}\}.
\end{equation*}
Equivalently, in the definition of $\JL(\mathcal{A})$, we restrict to finite cylinders rather than all finite subsets of $\omega \times 2$. 

It is worthwhile to briefly recall how $\JL(\mathcal{A})$ sits inside $\JL(\mathcal{B})$. Concretely, $\JL(\mathcal{A})$ consists of those functions in $\JL(\mathcal{B})$ which are constant on cylinders. The map $P: \JL(\mathcal{B}) \to \JL(\mathcal{B})$ defined by
\begin{equation*}
    Pf(n,0) = Pf(n,1) = \tfrac{1}{2}\big(f(n,0) + f(n,1)\big)
\end{equation*}
is a norm-one projection. We denote $Y = \ker P$, which can be explicitly described as
\begin{equation*}
    Y = \{f \in \JL(\mathcal{B}) : f(n,0) = -f(n,1) \text{ for all } n \in \omega\}.
\end{equation*}
Thus, we obtain the decomposition
\begin{equation}\label{eq: decomposition}
    C(K_\mathcal{B}) = \JL(\mathcal{B}) = \JL(\mathcal{A}) \oplus Y,
\end{equation}
which is precisely the content of \cite[Proposition~3.1]{plebanek}. Observe that by the discussion on Subsections \ref{subsec: JohnsonLind} and \ref{subsec: description-measures}, we can naturally identify
\begin{equation*}
    \JL(\mathcal{B})^* = M(\mathfrak{B}) =  C(K_\mathcal{B})^*  = \ell_1(\omega \times 2) \oplus \ell_1(\mathcal{B}) = \ell_1(\omega \times 2) \oplus \ell_1(\mathfrak{c} \times 2),
\end{equation*}
where $\mathfrak{B} = [\mathcal{B} \cup \finset{\omega \times 2}]$.

\subsection{Two spaces.}\label{subsec: two-spaces} As we mentioned before, we will need to carry out two simultaneous constructions and thus, to avoid confusion, it is convenient to use the distinctive notation $\tilde{\omega} \times 2$ for the second copy of $\omega \times 2$ in which the second construction will be carried out.

We emphasize that the absence of a tilde indicates we are working in $\omega \times 2$, whereas a tilde denotes objects in $\tilde{\omega} \times 2$. For instance, a cylinder $C$ is understood to satisfy $C \subseteq \omega \times 2$, while a cylinder $\tilde{C}$ satisfies $\tilde{C} \subseteq \tilde{\omega} \times 2$. It will also be convenient to sometimes use this when we are working with index sets. For example, we may write $(\tilde{c}_{\tilde{n}})_{\tilde{n} \in \tilde{\omega}}$. When there is no ambiguity, we may simply omit the tilde in the index. In all cases, the intended meaning should be clear from the context.

Therefore, in our construction, we will build two almost disjoint families $\mathcal{B} = \{B^0_\xi, B^1_\xi: \xi < \mathfrak{c} \} \subseteq \mathcal{P}(\omega \times 2)$ and $\tilde{\mathcal{B}} = \{\tilde{B}^0_{\xi}, \tilde{B}^1_{\xi}: \xi < \mathfrak{c} \} \subseteq \mathcal{P}(\tilde{\omega} \times 2)$ splitting almost disjoint families of cylinders $\mathcal{A} = \{A_\xi: \xi < \mathfrak{c}\}$ and $\tilde{\mathcal{A}} = \{\tilde{A}_{\xi}: \xi < \mathfrak{c}\}$. From the discussion above \eqref{eq: decomposition}, it follows that
\begin{equation*}
    \JL(\mathcal{B}) = \JL(\mathcal{A}) \oplus Y \hspace{5pt} \text{ and } \hspace{5pt} \JL(\tilde{\mathcal{B}}) = \JL(\tilde{\mathcal{A}}) \oplus \tilde{Y},
\end{equation*}
so that
\begin{equation*}
    \JL(\mathcal{B}) \oplus \JL(\tilde{\mathcal{B}}) =  (\JL(\mathcal{A}) \oplus Y) \oplus (\JL(\tilde{\mathcal{A}}) \oplus \tilde{Y}) = X \oplus \tilde{X},
\end{equation*}
where $X := \JL(\tilde{A}) \oplus Y$ and $\tilde{X} := \JL(A) \oplus \tilde{Y}$. Observe that this is precisely the setting in \eqref{eq: X-and-X}, bearing in mind that $\JL(\cdot)$ are isometric to $C(K)$-spaces.

For brevity, we write
\begin{equation*}
    \JL (\mathcal{B}, \tilde{\mathcal{B}}) = \JL (\mathcal{B}) \oplus \JL (\tilde{\mathcal{B}}) \hspace{5pt} \text{ and } \hspace{5pt} M(\mathfrak{B}, \tilde{\mathfrak{B}}) = \JL (\mathcal{B}, \tilde{\mathcal{B}})^* = M(\mathfrak{B}) \oplus M(\tilde{\mathfrak{B}}).
\end{equation*}
    We denote by $Q$ and $\tilde{Q}$ the norm-one projections of $\JL (\mathcal{B}, \tilde{\mathcal{B}})$ onto $X$ and $\tilde{X}$, respectively.

Furthermore, we denote by $M_1(\mathfrak{B}, \tilde{\mathfrak{B}})$ the unit ball of $\JL (\mathcal{B}, \tilde{\mathcal{B}})^*$, equipped with the weak$^*$ topology. We can naturally view $M(\mathfrak{B})$ and $M(\tilde{\mathfrak{B}})$ as subsets of $M(\mathfrak{B}, \tilde{\mathfrak{B}})$. In this spirit, any measure $\nu \in M(\mathfrak{B})$ can be regarded as an element of $M(\mathfrak{B}, \tilde{\mathfrak{B}})$, and similarly for $\tilde{\nu} \in M(\tilde{\mathfrak{B}})$. Naturally, any measure $\lambda \in M(\mathfrak{B}, \tilde{\mathfrak{B}})$ can be decomposed as $\lambda = \nu +\tilde{\nu}$ where $\nu \in M(\mathfrak{B})$ and $\tilde{\nu} \in M(\tilde{\mathfrak{B}})$.

\subsection{Norming sets and Banach lattices.}

Our overall aim is to prove that, for a suitable choice of the families $\mathcal{A}, \tilde{\mathcal{A}}, \mathcal{B}$ and $\tilde{\mathcal{B}}$, the resulting spaces $X$ and $\tilde{X}$ are not isomorphic to Banach lattices. Therefore, we recall the conditions from \cite{de2023counterexample} that ensure that this to occurs. We denote by $B_X$ the unit ball of a Banach space $X$.

\begin{definition}
    Let $X$ be a Banach space and $0 < c \leq 1$. 
    We say that a sequence $(e_n^*)_{n\in\omega} \subseteq B_{X^*}$ is $c$-norming if
    \begin{equation*}
        \norm{x} \geq c \sup_{n \in \omega} |\langle e^*_n, x\rangle|
    \end{equation*}
    for every $x \in X$. If the constant $c$ is not relevant, we simply say that the sequence $(e_n^*)_{n\in\omega} \subseteq B_{X^*}$ is norming.
\end{definition}

We have the following definition \cite[Definition 4.2]{de2023counterexample}, which encapsulates the fact that the absolute value of $f$ does not exists.

\begin{definition}
    We say that a Banach space $X$ has the \emph{Desired Property of De Hevia, Mart\'inez-Cervantes, Salguero-Alarc\'on, and Tradacete}, abbreviated (DP), if for every norming sequence 
    $(e_n^*)_{n\in\omega}$ in $B_{X^*}$ there exists an element $f \in X$ such that no element $g \in X$ satisfies
    \begin{equation}\label{eq: DP}
       \langle e_n^*, g \rangle = |\langle e_n^*, f \rangle| \quad \text{for every } n \in \omega. \tag{DP}
    \end{equation}
\end{definition}

The motivation behind the previous definition is that, under certain conditions, satisfying \eqref{eq: DP} implies that the space is not isomorphic to a Banach lattice. In particular, we have the following, see \cite[Corollary 4.3]{de2023counterexample}.

\begin{proposition}\label{lmm: non-banach-lattice-criteria}
    Let $X$ be a Banach space with a countable norming set. Suppose that $X^*$ is isomorphic to $\ell_{1}(\Gamma)$ for some index set $\Gamma$. Then, the following statements are equivalent:
    \begin{enumerate}
        \item $X$ is isomorphic to a Banach lattice.
        \item $X$ is isomorphic to a sublattice of $\ell_{\infty}$.
        \item $X$ does not have the \eqref{eq: DP}. That is, there exists a norming sequence $(e_n^*)_{n\in\omega}$ in $B_{X^*}$ such that for every $f \in X$ there is an element $g \in X$ such that
        \begin{equation*}
           \langle e_n^*, g \rangle = |\langle e_n^*, f \rangle|, \quad \text{for every } n \in \omega.
        \end{equation*}
    \end{enumerate}
\end{proposition}

We note that, for any almost disjoint families $\mathcal{A}, \tilde{\mathcal{A}}, \mathcal{B}$ and $\tilde{\mathcal{B}}$, the spaces $X$ and $\tilde{X}$, defined in Subsection \ref{subsec: two-spaces}, satisfy the conditions of Proposition \ref{lmm: non-banach-lattice-criteria} (see Section \ref{sec: lattice} for details). Therefore, to prove they are not isomorphic to Banach lattices, it is enough to show that they satisfy \eqref{eq: DP}.

\subsection{The spaces $M(\mathfrak{B}, \tilde{\mathfrak{B}})$, $X^*$ and $\tilde{X}^*$.}\label{subsec: measures}

Let $\mathcal{A}, \tilde{\mathcal{A}}, \mathcal{B}$ and $\tilde{\mathcal{B}}$ be almost disjoint families and $\JL(\mathcal{B}, \mathcal{B}) = X \oplus \tilde{X}$ as in the setting of Subsection \ref{subsec: two-spaces}. To ensure that $X$ and $\tilde{X}$ are not isomorphic to Banach lattices, it will be necessary to have control over the elements in the duals $X^*$ and $\tilde{X}^*$. 

We start with the following observation, which is the analogue of \cite[Remark~3.2]{plebanek}. Recall that $Q$ and $\tilde{Q}$ denote the projections onto $X$ and $\tilde{X}$, respectively.

\begin{rem}\label{rem: norming-as-subsets-measures}
    Let $x^* \in X^*$. Then $x^* Q \in \JL (\mathcal{B}, \tilde{\mathcal{B}})^* = M(\mathfrak{B}, \tilde{\mathfrak{B}})$, so that every functional on $X$ can be seen as an element on $M(\mathfrak{B}, \tilde{\mathfrak{B}})$ vanishing on all cylinders from $\mathfrak{B}$. In other words, we have an embedding
    \begin{equation*}
        X^* \to M(\mathfrak{B}, \tilde{\mathfrak{B}}), \quad x^* \mapsto x^* Q.
    \end{equation*}
    Analogous considerations apply to $\tilde{x}^* \in \tilde{X}^*$. 
\end{rem}

Thus, from now on, when we consider an element in the dual of $X$ or $\tilde{X}$, we will always think of it as the measure in $M(\mathfrak{B}, \tilde{\mathfrak{B}})$ representing it. 

\begin{rem}\label{rem: decomp-of-measures}
    We will reserve the symbols $\lambda$ and $\gamma$ for elements in $M_1(\mathfrak{B}, \tilde{\mathfrak{B}})$. Observe that we have a decomposition $\lambda = \nu + \tilde{\nu}$ where $\nu \in M_1(\mathfrak{B})$ and $\tilde{\nu} \in M_1(\tilde{\mathfrak{B}})$ (with their natural interpretations).  From the discussion in Subsection \ref{subsec: description-measures}, we can further express $\nu = \mu + \overline{\nu}$ where $\mu \in \ell_1(\omega \times 2)$, $\overline{\nu}$ vanishes on finite sets and $|\overline{\nu}|(B) \not = 0$ for at most countably many $B \in \mathcal{B}$. Identical considerations apply for $\tilde{\nu}$. Using these decompositions we can write $\lambda = \nu + \tilde{\nu} = (\mu + \overline{\nu}) + (\tilde{\mu} + \tilde{\overline{\nu}})$. Therefore, a measure $\lambda \in M_1(\mathfrak{B}, \tilde{\mathfrak{B}})$ can be seem as a quadruple $(\mu, \overline{\nu}; \tilde{\mu}, \tilde{\overline{\nu}})$. We call $\mu$ and $\tilde{\mu}$ the $\ell_1(\omega\times2)$-part and $\ell_1(\tilde{\omega} \times 2)$-part of $\lambda$ respectively. 

    Furthermore, from the discussion in Subsection~\ref{subsec: description-measures}, 
    we can identify $\mu$ and $\tilde{\mu}$ with elements 
    $x \in \ell_1(\omega \times 2)$ and $\tilde{x} \in \ell_1(\tilde{\omega} \times 2)$, 
    and $\overline{\nu}$ and $\tilde{\overline{\nu}}$ with elements 
    $y \in \ell_1(\mathfrak{c} \times 2)$ and $\tilde{y} \in \ell_1(\tilde{\mathfrak{c}} \times 2)$. 
    In this way, the measure $\lambda$ can be coded by the quadruple
    \begin{equation*}
       (x,y;\,\tilde{x},\tilde{y}) \in 
       \ell_1(\omega \times 2) \oplus \ell_1(\mathfrak{c} \times 2) 
       \oplus \ell_1(\tilde{\omega} \times 2) \oplus \ell_1(\tilde{\mathfrak{c}} \times 2).
    \end{equation*}
\end{rem}

\bigskip
\section{Technical results}\label{sec: technical-results}

We recall some results of Plebanek and Salguero-Alarc\'on and make the natural adjustments to fit our framework.

\subsection{Separating sets of measures.}

Since the nature of the underlying countable set is not important, in this section, we will discuss measures on $\omega$ and $\tilde{\omega}$. Throughout this section, we fix Boolean subalgebras 
$\mathfrak{B} \subseteq \mathcal{P}(\omega)$ and 
$\tilde{\mathfrak{B}} \subseteq \mathcal{P}(\tilde{\omega})$ 
containing all finite subsets of $\omega$ and $\tilde{\omega}$, respectively. We state our results in general form, for arbitrary finitely-additive measures $\lambda$ on $\omega \cup \tilde{\omega}$, which we implicitly regard as the sum of a measure on $\omega$ and a measure on $\tilde{\omega}$. Accordingly, we adopt the slight abuse of notation $M_1(\mathcal{P}(\omega), \mathcal{P}(\tilde{\omega}))$ for the collection of those measures with total variation at most one.

We shall employ the following analogue of separation, as given in \cite[Definition 4.1]{plebanek}.

\begin{definition}
    Let $M, M' \subseteq M_1(\mathcal{P}(\omega), \mathcal{P}(\tilde{\omega}))$. We say that $(M, M')$ is $(\mathfrak{B}, \tilde{\mathfrak{B}})$-separated with constant $\varepsilon$ if there exist $\varepsilon > 0$, $n, m \in \omega$ and $B_1, \dots, B_n \in \mathfrak{B}$, $\tilde{B}_1, \dots, \tilde{B}_m \in \tilde{\mathfrak{B}}$ such that
    \begin{equation*}
        \max( \max_{i \leq n} |\lambda(B_i) - \lambda'(B_i)|, \max_{j \leq m} |\lambda(\tilde{B}_j) - \lambda'(\tilde{B}_j)|) \geq \varepsilon
    \end{equation*}
    for any $\lambda \in M$ and $\lambda' \in M'$. If the constant is not relevant, we simply say that $M, M'$ are $(\mathfrak{B}, \tilde{\mathfrak{B}})$-separated.
\end{definition}

With this definition, we can get the following analogue of \cite[Lemma 4.2]{plebanek}. The proof is identical and thus omitted.

\begin{lemma}\label{lmm: simple-functions-separate-sets}
    Let $M, M' \subseteq  M_1(\mathcal{P}(\omega), \mathcal{P}(\tilde{\omega}))$ and suppose there exist $\varepsilon > 0$ and simple $\mathfrak{B}$-measurable and $\tilde{\mathfrak{B}}$-measurable functions $g: \omega \to \mathbb{R}$ and $\tilde{g}: \tilde{\omega} \to \mathbb{R}$ such that
    \begin{equation*}
        |\langle \lambda, g + \tilde{g}\rangle - \langle \lambda', g + \tilde{g}\rangle| \geq \varepsilon,
    \end{equation*}
    for any $\lambda \in M$ and $\lambda' \in M'$. Then $(M,M')$ is $(\mathfrak{B}, \tilde{\mathfrak{B}})$-separated.
\end{lemma}

Likewise, we have the following counterpart of \cite[Lemma 4.3]{plebanek}. The proof is essentially the same, using now the compactness of $[-1,1]^n \times [-1,1]^m$. We omit the details.

\begin{lemma} \label{lmm: only-countable-many}
    Suppose that $|\mathfrak{B}| < \mathfrak{c}$ and $|\tilde{\mathfrak{B}}| < \mathfrak{c}$. Then, for every infinite subset $M \subseteq M_1(\mathcal{P}(\omega), \mathcal{P}(\tilde{\omega}))$ the pair $(M', M\setminus M')$ can be $(\mathfrak{B}, \tilde{\mathfrak{B}})$-separated for fewer than $\mathfrak{c}$ many sets $M' \subseteq M$.
\end{lemma}

For any set $Z \subseteq \omega$, we write $\mathfrak{B}[Z]$ for the algebra generated by $\mathfrak{B} \cup \{Z\}$. Observe that every element in $\mathfrak{B}[Z]$ is of the form $(A \cap Z) \cup (B \cap Z^c)$ for some $A, B \in \mathfrak{B}$.  Similar for $\tilde{Z} \subseteq \tilde{\omega}$ and $\tilde{\mathfrak{B}}$.

Finally, we will need the following counterpart of \cite[Proposition 4.4]{plebanek}, whose underlying ideas date back to a counting argument of Haydon~\cite{haydon1981non}. 

\begin{proposition}\label{prop: Haydon}
    Suppose $I$ is an index set with $|I| < \mathfrak{c}$ and that we are given:
    \begin{enumerate}[label=(\roman*)]
        \item a list $\{(\mathfrak{B}_j, \tilde{\mathfrak{B}}_j): j \in I\}$ of pairs of Boolean subalgebras of $\mathcal{P}(\omega)$ and $\mathcal{P}(\tilde{\omega})$ respectively, such that, for every $j \in I$, $\finset{\omega} \subseteq \mathfrak{B}_j$, $\finset{\tilde{\omega}} \subseteq \tilde{\mathfrak{B}}_j$ and $|\mathfrak{B}_j| \leq |I| $, $|\tilde{\mathfrak{B}}_j| \leq |I|$.
        \item \label{it: part2} a list $\{(M_j, M'_j): j \in I \}$ of pairs of subsets of $M_1(\mathcal{P}(\omega), \mathcal{P}(\tilde{\omega}))$ such that, for every $j \in J$, the pair $(M_j, M'_j)$ is not $(\mathfrak{B}_j, \tilde{\mathfrak{B}}_j)$-separated.
    \end{enumerate}
    Then for every almost disjoint family $\mathcal{Z} \subseteq \mathcal{P(\omega)}$ with $|\mathcal{Z}| = \mathfrak{c}$ there exists $Z \in \mathcal{Z}$ such that for every $j \in J$, the pair $(M_j, M'_j)$ is not $(\mathfrak{B}_j[Z], \tilde{\mathfrak{B}}_j)$-separated. An analogous result holds for any almost disjoint family $\tilde{\mathcal{Z}} \subseteq \mathcal{P(\tilde{\omega})}$ with $|\tilde{\mathcal{Z}}| = \mathfrak{c}$.
\end{proposition}

We give a brief outline of the proof, which is essentially the same as the proof of \cite[Proposition 4.4]{plebanek}. We do this for the almost disjoint family $\mathcal{Z}$, the other case being identical. We start with an analogue to \cite[Lemma 4.5]{plebanek}.
\begin{lemma}\label{lmm: 4-5counter}
    Suppose that for a given $M, M' \subseteq M_1(\mathcal{P}(\omega), \mathcal{P}(\tilde{\omega}))$, there are:
    \begin{enumerate}
        \item $n \in \omega$, $m \in \omega$, $\varepsilon > 0$ and $k > \frac{12n}{\varepsilon}$
        \item almost disjoint sets $Z_1, \dots, Z_k \subseteq \omega$,
        \item $A_1, B_1, \dots, A_n, B_n \in \mathfrak{B}$ and $\tilde{A}_1, \dots, \tilde{A}_m \in \tilde{\mathfrak{B}}$
    \end{enumerate}
     such that, for every $1 \leq j \leq k$, the pair $(M, M')$ is $(\mathfrak{B}[Z_j], \tilde{\mathfrak{B}})$-separated, with constant $\varepsilon$, by the sets
     \begin{equation*}
         Y_{j,i} = (A_i \cap Z_j) \cup (B_i \cap Z_j^c), \hspace{5pt} 1 \leq i \leq n \hspace{10 pt} \text{ and } \tilde{A}_l, \hspace{10pt} 1 \leq l \leq m.
     \end{equation*}
     Then $(M, M')$ is $(\mathfrak{B}, \tilde{\mathfrak{B}})$-separated.
\end{lemma}
\begin{proof}[Outline of the proof]
    Consider a fixed pair $(\lambda, \lambda') \in M \times M'$. If there exists $1 \leq l \leq m$ such that $|\lambda(\tilde{A}_l) - \lambda'(\tilde{A}_l) | \geq \varepsilon$, then we are done. Otherwise, we can argue word for word as in the proof of \cite[Lemma 4.5]{plebanek}.
\end{proof}

Once we have this analogue of \cite[Lemma 4.5]{plebanek}, we can prove the extension to $(\mathfrak{B}_j[Z], \tilde{\mathfrak{B}}_j)$ in Proposition \ref{prop: Haydon}, by Haydon's counting argument. Suppose that for every $Z \in \mathcal{Z}$ the conclusion of Proposition \ref{prop: Haydon} fails. This is witnessed by some $j \in I$, $n \in \omega$, $m \in \omega$, some rational $\varepsilon > 0$, and some collections $A_1, B_1, \dots, A_n, B_n \in \mathfrak{B}_j$ and $\tilde{A}_1, \dots, \tilde{A}_m \in \tilde{\mathfrak{B}}_j$. Since the number of such witnesses has size less than $\mathfrak{c}$ and $|\mathcal{Z}| = \mathfrak{c}$, there are infinitely many $Z \in \mathcal{Z}$ with the same witness. Lemma \ref{lmm: 4-5counter} gives that the pair $(M_j, M'_j)$ is $(\mathfrak{B}_j, \tilde{\mathfrak{B}}_j)$-separated, contrary to our assumptions.

\subsection{Separation versus norming sets.} By definition, the space $X$ contains every function of the form $f_n =  \mathds{1}_{(n,0)} - \mathds{1}_{(n,1)} \in Y \subseteq X$ for $n \in \omega$, while the space $\tilde{X}$ contains every function of of the form $\tilde{f}_{\tilde{n}} = \mathds{1}_{(\tilde{n},0)} - \mathds{1}_{(\tilde{n},1)} \in \tilde{Y} \subseteq \tilde{X}$ for $\tilde{n} \in \tilde{\omega}$. Regarding elements in $B_{X^*}$ and $B_{\tilde{X}^*}$ as elements in $M_1(\mathfrak{B}, \tilde{\mathfrak{B}})$, as in Remark \ref{rem: norming-as-subsets-measures}, then the elements in $B_{X^*}$ vanish on all the cylinders from $\mathfrak{B}$ while the elements in $B_{\tilde{X}^*}$ vanish on all cylinders from $\tilde{\mathfrak{B}}$.

From the above discussion, every norming sequence for $X$, regarded as a sequence in $M_1(\mathfrak{B}, \tilde{\mathfrak{B}})$, must contain a subsequence $(\lambda_n)_{n \in \omega}$ such that $\inf_{n \in \omega}|\lambda_n(\{(n,0)\})| = \frac{1}{2}\inf_{n \in \omega}|\langle \lambda_n, f_n \rangle|  > 0$. Observe that if we decompose $\lambda = (\mu_n + \overline{\nu}_n) + \tilde{\nu}_n$ as in Remark \ref{rem: decomp-of-measures}, then $\lambda_n(\{(n,0)\}) = \mu_n(\{(n, 0)\})$. Similar considerations apply for $\tilde{X}$; this motivates the following, see \cite[Definition 5.1]{plebanek}.

\begin{definition}\label{def: admissible-sequences}
    We say that a sequence $(\mu_n)_{n\in \omega}$ in the unit ball of $\ell_1(\omega \times 2)$ is \emph{admissible} if 
    \begin{enumerate}
        \item $\mu_n(\{(k,0)\}) = - \mu_n(\{(k,1)\})$ for all $n \in \omega$ and $k \in \omega$,
        \item $\inf_{n \in \omega} |\mu_n(\{(n, 0)\})| > 0$.
    \end{enumerate}
    Similarly, a sequence $(\tilde{\mu}_{\tilde{n}})_{\tilde{n} \in \tilde{\omega}}$ in the unit ball of $\ell_1(\tilde{\omega} \times 2)$ is $\sim$-\emph{admissible} if 
    \begin{enumerate}
        \item $\tilde{\mu}_{\tilde{n}}(\{(\tilde{k},0)\}) = - \tilde{\mu}_{\tilde{n}}(\{(\tilde{k},1)\})$ for all $\tilde{n} \in \tilde{\omega}$ and $\tilde{k} \in \tilde{\omega}$,
        \item $\inf_{\tilde{n} \in \tilde{\omega}} |\tilde{\mu}_{\tilde{n}}(\{(\tilde{n}, 0)\})| > 0$.
    \end{enumerate}
\end{definition}

For a sequence $(\mu_n)_{n \in \omega}$ and any set $J \subseteq \omega$ we write $\mu[J] = \{ \mu_n: n \in J\}$. For convenience, we introduce the following notation: for $s, t \in \mathbb{R}$ and $\varepsilon > 0$ we say $s \approx_\varepsilon t$ if $|s - t| < \varepsilon$. We have the following version of \cite[Lemma 5.3]{plebanek}. 

\begin{lemma}\label{lmm: main}
    Let $C = C_0 \times 2 \subseteq \omega \times 2$ 
    be a cylinder, $(\mu_n)_{n \in \omega}$ be an admissible sequence and $b = \inf_{n \in \omega} |\mu_n(\{(n,0)\})|$. Then there exist:
    \begin{itemize}[label = --]
        \item three infinite sets $J_2 \subseteq J_1 \subseteq J_0 \subseteq C_0$ with $J_0 \setminus J_1$, $J_1 \setminus J_2$ infinite,
        \item a set $Z$ splitting the cylinder $J_1 \times 2$,
        \item two positive constants $a \geq b$ and $0 < \delta < b/11$,
    \end{itemize}
    such that the following holds:
    \begin{equation}\label{eq: A1}
        \mu_n(Z) - \mu_n((J_1 \times 2 )\setminus Z) \approx_{3\delta} \begin{cases}
            2a & \text{if } n \in J_2,  \\
            -2a & \text{if } n \in J_1 \setminus J_2, \\
            0 & \text{if } n \in J_0 \setminus J_1.
        \end{cases} \tag{A1}
    \end{equation}
    Similar considerations apply for a cylinder $\tilde{C} = \tilde{C}_0 \times 2 \subseteq \tilde{\omega} \times 2$ and a $\sim$-admissible sequence $(\tilde{\mu}_{\tilde{n}})_{\tilde{n} \in \tilde{\omega}}$.
\end{lemma}
\begin{proof}
    The proof is identical to that of \cite[Lemma 5.3]{plebanek}. Since it is short and crucial for our purposes, we include it for the reader's convenience.

    Choose $0 < \delta < b/11$. By Rosenthal's lemma (see \cite[p. 18]{rosenthal1970relatively} and the version \cite[Theorem 3.5]{plebanek}), we can find an infinite subset $J_0 \subseteq C_0$ such that
    \begin{equation*}
        |\mu_n|((J_0 \times 2) \setminus c_n) < \delta \hspace{5 pt} \text{ for every } n \in J_0.
    \end{equation*}
    We can now find an infinite subset $J_1 \subseteq J_0$ with $J_0 \setminus J_1$ infinite, and singletons $p_n \subseteq c_n$ so that the sequence $(\mu_n(p_n))_{n \in J_1}$ converges to some $a \geq b$. Observe that these are precisely the singletons such that $\mu_n(p_n) > 0$.
    
    By removing a finite number of elements in $J_1$, we may assume that
    \begin{equation*}
        |\mu_n(p_n) - a| < \delta  \text{ for every } n \in J_1.
    \end{equation*}
    Lastly, choose an infinite set $J_2 \subseteq J_1$ with $J_1 \setminus J_2$ infinite and set
    \begin{equation*}
        Z = \left(\bigcup_{n \in J_2} p_n \right) \cup \left( \bigcup_{n \in J_1 \setminus J_2} (c_n \setminus p_n)\right). 
    \end{equation*}
    By construction, $Z$ splits the cylinder $J_1 \times 2$, so that we only need to verify that \eqref{eq: A1} is satisfied. For every $n \in J_0 \backslash J_1$ we have $J_1 \times 2 \subseteq (J_0 \times 2) \setminus c_n$, which gives
    \begin{align*}
        |\mu_n(Z) - \mu_n((J_1 \times 2) \setminus Z)| \leq |\mu_n|((J_0 \times 2) \setminus c_n) < \delta,
    \end{align*}
    equivalently $\mu_n(Z) - \mu_n((J_1 \times 2) \setminus Z) \approx_\delta 0$.

    Observe that for any $n \in J_1$, we have $\mu_n(p_n) \approx_\delta a$ and $\mu_n(c_n \setminus p_n) \approx_\delta -a$. Now, for every $n \in J_1 \backslash J_2$, we have $(c_n \backslash p_n) \in Z$ and $p_n \in (J_1 \times 2) \setminus Z$ so that
    \begin{equation*}
        \mu_n(Z) - \mu_n((J_1 \times 2) \setminus Z) \approx_\delta \mu_n(c_n \backslash p_n)  - \mu_n(p_n)  \approx_{2\delta} -2a.
    \end{equation*}
    Similarly, for $n \in J_2$ we have $p_n \in Z$ and $(c_n \setminus p_n) \in (J_1 \times 2) \setminus Z$ and thus
    \begin{equation*}
        \mu_n(Z) - \mu_n((J_1 \times 2) \setminus Z) \approx_\delta \mu_n( p_n)  - \mu_n( c_n \setminus p_n)  \approx_{2\delta} 2a,
    \end{equation*}
    which finishes the proof.
\end{proof}

We also recall the following, which is \cite[Remark 5.4]{plebanek}, and follows by examining the proof of the previous lemma.

\begin{rem}\label{rem: uncountably-many-ways}
    In Lemma \ref{lmm: main}, the sets $J_1$ and $J_2$ may be replaced by any infinite subsets of themselves; in particular, there are continuum many ways to choose them. Observe that, if we do this replacement, then \eqref{eq: A1} is still satisfied with the same values of $\delta$ and $a$.
\end{rem}

\bigskip
\section{The construction} \label{sec: proof-main-theorem}

\subsection{Notation and Setup.} We now present the construction of the four almost disjoint families $\mathcal{A}$, $\tilde{\mathcal{A}}$, $\mathcal{B}$ and $\tilde{\mathcal{B}}$, having the properties specified in Subsection \ref{subsec: two-spaces}. We closely follow the construction in \cite[Section 6]{plebanek}, with modifications to adapt it to our setting as appropriate.

We will construct almost disjoint families $\mathcal{A} = \{A_\alpha: \alpha < \mathfrak{c} \}$ and $\tilde{\mathcal{A}} = \{\tilde{A}_\alpha: \alpha < \mathfrak{c} \}$ of cylinders of $\omega \times 2$ and $\tilde{\omega} \times 2$, and splittings of every $A_\alpha$ and $\tilde{A}_\alpha$ into two sets $B_\alpha^{i}$ and $\tilde{B}_\alpha^{i}$, $i=0,1$. The construction will be carried out so that $X = \JL(\tilde{\mathcal{A}}) \oplus Y$ and $\tilde{X} = \JL(\mathcal{A}) \oplus \tilde{Y}$ fail to be isomorphic to Banach lattices.

We start by coding all the sequences of measures that will (eventually) appear in norming sequences of $X^*$ and $X^*$. For $(x,y; \tilde{x}, \tilde{y}) \in \ell_1(\omega \times 2) \times \ell_1(\mathfrak{c} \times 2) \times \ell_1(\tilde{\omega} \times 2) \times \ell_1(\tilde{\mathfrak{c}} \times 2)$ we say that $(x,y; \tilde{x}, \tilde{y})$ is of
\begin{enumerate}[label = Type-\arabic*, ref =  Type-\arabic*,leftmargin=*]
        \item\label{type-1} if $x(n,0) + x(n,1) = 0$ for every $n \in \omega$ and $y(\alpha, 0) + y(\alpha,1) = 0$ for every $\alpha < \mathfrak{c}$,
        \item\label{type-2} if $\tilde{x}(\tilde{n},0) + \tilde{x}(\tilde{n},1) = 0$ for every $\tilde{n} \in \tilde{\omega}$ and $\tilde{y}(\tilde{\alpha}, 0) + \tilde{y}(\tilde{\alpha},1) = 0$ for every $\tilde{\alpha} < \tilde{\mathfrak{c}}$.
\end{enumerate}

We let $\Psi$ be the  set of all $(x,y; \tilde{x}, \tilde{y})$ of \ref{type-1} such that additionally $\norm{x}_1 + \norm{y}_1 + \norm{\tilde{x}}_1 + \norm{\tilde{y}}_1 \leq 1$. Similarly $\tilde{\Psi}$ is the set of all $(x,y; \tilde{x}, \tilde{y})$ of \ref{type-2} satisfying $\norm{x}_1 + \norm{y}_1 + \norm{\tilde{x}}_1 + \norm{\tilde{y}}_1 \leq 1$. We will only care about those sequences of measures that (will eventually) appear as norming sequences of $X^*$ and $\tilde{X}^*$, and thus we will make a further reduction. 

Hence we let $\Phi$ be the set of all sequences $(x_n,y_n; \tilde{x}_n, \tilde{y}_n)_{n \in \omega}$ in $\Psi$ such that
\begin{equation*}
    \inf_{n \in \omega} |x_n(n,0)| > 0,
\end{equation*}
and similarly we let $\tilde{\Phi}$ be the set of all sequences $(x_{\tilde{n}},y_{\tilde{n}}; \tilde{x}_{\tilde{n}}, \tilde{y}_{\tilde{n}})_{{\tilde{n}} \in \tilde{\omega}}$ in $\tilde{\Psi}$ such that
\begin{equation*}
    \ \inf_{\tilde{n} \in \tilde{\omega}} |\tilde{x}_{\tilde{n}}(\tilde{n},0)| > 0. \\
\end{equation*}

We can enumerate 
\begin{equation*}
    \Phi = \{\varphi_\alpha : \alpha < \mathfrak{c}\} \quad 
\end{equation*}
in such a way that, if $\varphi_\alpha = (x^\alpha_n, y^\alpha_n; \tilde{x}^\alpha_n, \tilde{y}^\alpha_n)_{n \in \omega} \in \Phi$, then
\begin{equation}\label{eq: enumeration}
    y^\alpha_n(\xi, i) = \tilde{y}^\alpha_n(\tilde{\xi}, i) = 0 \quad \text{for every } n \in \omega, i = 0,1 \text{ and every } \xi \geq \alpha \text{ and } \ \tilde{\xi} \geq \tilde{\alpha}.
\end{equation}

To achieve this, we follow the argument in \cite[Section 6]{plebanek}. Start with an enumeration $\{\theta_\alpha: \alpha < \mathfrak{c} \}$ of $\Phi$ where every element is repeated cofinally often. Then for $\theta_\alpha = (x^\alpha_n, z^\alpha_n; \tilde{x}^\alpha_n, \tilde{z}^\alpha_n)_{n \in \omega}$ put $\varphi_\alpha = (x^\alpha_n,y^\alpha_n; \tilde{x}^\alpha_n, \tilde{y}^\alpha_n)_{n \in \omega}$ where for every $n \in \omega$ and $i= 0,1$, $y_n^\alpha(\beta, i) = z^\alpha_n(\beta,i)$ and $\tilde{y}^\alpha_n(\tilde{\beta},i) = \tilde{z}^\alpha_n(\tilde{\beta},i)$ whenever $\beta < \alpha$ and $\tilde{\beta} < \tilde{\alpha}$ and $y^\alpha_n(\beta,i) = \tilde{y}^\alpha_n(\tilde{\beta},i) = 0$ whenever $\beta \geq \alpha$ and $\tilde{\beta} \geq \tilde{\alpha}$.

Similarly, we can enumerate
\begin{equation*}
    \tilde{\Phi} = \{\tilde{\varphi}_{\alpha} : \alpha < \mathfrak{c}\}
\end{equation*} 
such that if $\tilde{\varphi}_\alpha = (x^\alpha_{\tilde{n}}, y^\alpha_{\tilde{n}}; \tilde{x}^\alpha_{\tilde{n}}, \tilde{y}^\alpha_{\tilde{n}})_{{\tilde{n}} \in \tilde{\omega}} \in \Phi$ then \eqref{eq: enumeration} holds for every $\tilde{n} \in \tilde{\omega}$. \\

We can now exploit the following devilishly clever observation of Plebanek and Salguero-Alarc\'on: using these enumerations, the codings $\varphi_\alpha$ and $\varphi'_\alpha$ represent sequences of measures on $M_1(\mathfrak{B}, \tilde{\mathfrak{B}})$ \emph{even before} the construction of $\mathcal{B}$ and $\mathcal{B}'$ has been completed. This is a delicate point which warrants discussion. 

At induction step $\alpha$, the sets $\{B_\beta^0, B_\beta^1: \beta < \alpha\}$ and $\{\tilde{B}_\beta^0, \tilde{B}_\beta^1: \beta < \alpha\}$ will have been constructed, and thus we can consider the Boolean algebras that they generate
\begin{equation*}
    \mathfrak{B}(\alpha) = [\{B_{\beta}^0, B_{\beta}^1: \beta < \alpha  \} \cup \finset{\omega \times 2}] \hspace{5pt} \text{ and }  \tilde{\mathfrak{B}}(\alpha) = [\{\tilde{B}_{\beta}^0, \tilde{B}_{\beta}^1: \beta < \alpha \} \cup \finset{\omega \times 2}].
\end{equation*}
Observe now that the element $\varphi_\alpha = (x^\alpha_n, y^\alpha_n; \tilde{x}^\alpha_n, \tilde{y}^\alpha_n)_{n \in \omega}$ can, at this point, be considered to code a sequence $(\lambda_n^\alpha)_{n \in \omega}$ of measures in $M_1(\mathfrak{B}(\alpha), \tilde{\mathfrak{B}}(\alpha))$ since \eqref{eq: enumeration} holds. Finally, the coup de grâce: \emph{once the construction} of $\mathcal{B}$ and $\tilde{\mathcal{B}}$ is complete, the algebras $\mathfrak{B}(\alpha)$ and $\tilde{\mathfrak{B}}(\alpha)$ will be Boolean subalgebras of $\mathfrak{B}$ and $\tilde{\mathfrak{B}}$. Thus, arguing as in Subsection \ref{subsec: description-measures}, the sequence of measures $(\lambda_n^\alpha)_{n \in \omega}$ in $M_1(\mathfrak{B}(\alpha), \tilde{\mathfrak{B}}(\alpha))$ extends uniquely to a sequence of measures in $M_1(\mathfrak{B}, \tilde{\mathfrak{B}})$, which, in an abuse of notation, we still call $(\lambda_n^\alpha)_{n \in \omega}$. This is precisely what we mean when we say that $\varphi_\alpha$ codes a sequence of measure $(\lambda_n^\alpha)_{n \in \omega}$ in $M_1(\mathfrak{B}, \tilde{\mathfrak{B}})$.

Conversely, observe that any sequence in $M_1(\mathfrak{B}, \tilde{\mathfrak{B}})$ arises as the extension of a sequence in $M_1(\mathfrak{B}(\alpha), \tilde{\mathfrak{B}}(\alpha))$ for some $\alpha < \mathfrak{c}$. For details of how this works, we refer back to Subsection \ref{subsec: description-measures}, together with the observation that the cofinality of the continuum cannot be $\omega$.

During the inductive construction, for a subset $V \subseteq \mathfrak{c}$ we will denote the Boolean subalgebras
\begin{equation*}
    \mathfrak{B}(V) = [\{B_{\beta}^0, B_{\beta}^1: \beta \in V  \} \cup \finset{\omega \times 2}] \hspace{5pt} \text{ and }  \tilde{\mathfrak{B}}(V) = [\{\tilde{B}_{\beta}^0, \tilde{B}_{\beta}^1: \beta \in V \} \cup \finset{\omega \times 2}],
\end{equation*}
provided all the elements $B_{\beta}^i, \tilde{B}_{\beta}^i$, $i = 0,1$, $\beta \in V$, have been constructed at that point of the induction.

\subsection{The main construction.} 

Before we start, we fix two almost disjoint families 
$\mathcal{R} = \{R_\alpha : \alpha < \mathfrak{c} \}$ and $\tilde{\mathcal{R}} = \{\tilde{R}_\alpha : \alpha < \mathfrak{c} \}$ of infinite subsets of $\omega$ and $\tilde{\omega}$ respectively.

From now on, we will write $(\lambda^\alpha_{n})_{n \in \omega}$ and $(\gamma^\alpha_{\tilde{n}})_{\tilde{n} \in \tilde{\omega}}$ to refer to the sequences in $M_1(\mathfrak{B}, \tilde{\mathfrak{B}})$ coded by $\varphi_\alpha$ and $\tilde{\varphi}_{\alpha}$ respectively. This should be understood in the sense explained above. For these sequences, we write
    \begin{equation*}
           \lambda^\alpha_{n} = (\mu^\alpha_n + \overline{\nu}^\alpha_n) + \tilde{\nu}^\alpha_n \hspace{5pt} \text{ and } \hspace{5pt} \gamma^\alpha_{\tilde{n}}= \nu^\alpha_{\tilde{n}} + (\tilde{\mu}^\alpha_{\tilde{n}} + \tilde{\overline{\nu}}^\alpha_{\tilde{n}}),
    \end{equation*}
as in Remark \ref{rem: decomp-of-measures}. \\

Using transfinite induction, we will construct for every $\xi < \mathfrak{c}$:
\begin{enumerate}[label=(\roman*), ref = (\roman*)]
        \item\label{it: enu1} three infinite subsets $J_2^\xi \subseteq J_1^\xi \subseteq J_0^\xi \subseteq R_\xi$,
        \item\label{it: enu2} three infinite subsets $\tilde{J}_2^\xi \subseteq \tilde{J}_1^\xi \subseteq \tilde{J}_0^\xi \subseteq \tilde{R}_\xi$,
        \item\label{it: enu3} splittings $B_\xi^0$, $B_\xi^1$ and $\tilde{B}_\xi^0$, $\tilde{B}_\xi^1$ of the cylinders  $A_\xi =J_1^\xi \times 2$ and $\tilde{A}_\xi = \tilde{J}_1^\xi \times 2$,
        \item\label{it: enu5} positive constants $a_\xi \geq b_\xi := \inf_{n \in \omega} |\mu^\xi_n(\{(n,0)\})|$, $\tilde{a}_\xi \geq \tilde{b}_\xi := \inf_{\tilde{n} \in \tilde{\omega}} |\tilde{\mu}^\xi_{\tilde{n} }(\{(\tilde{n} ,0)\})|$, $0 < \delta_\xi < b_\xi/11$ and $0 < \tilde{\delta}_\xi < \tilde{b}_\xi/11$;
\end{enumerate}
    so that the following conditions hold: \\
 \begin{enumerate}[label = (B\arabic*), ref = (B\arabic*)]
     \item \label{it: B1} \begin{equation*}
        \mu^\xi_n(B_\xi^0) - \mu^\xi_n(B_\xi^1) \approx_{3\delta_\xi} \begin{cases}
            2a_\xi & \text{if } n \in J^\xi_2,  \\
            -2a_\xi & \text{if } n \in J^\xi_1 \setminus J^\xi_2, \\
            0 & \text{if } n \in J^\xi_0 \setminus J^\xi_1.
        \end{cases}
    \end{equation*}
     \item \label{it: B2} \begin{equation*}
        \tilde{\mu}^\xi_{\tilde{n}}(\tilde{B}_\xi^0) - \tilde{\mu}^\xi_{\tilde{n}}(\tilde{B}_\xi^1) \approx_{3\tilde{\delta}_\xi} \begin{cases}
            2\tilde{a}_\xi & \text{if } \tilde{n} \in \tilde{J}^\xi_2,  \\
            -2\tilde{a}_\xi & \text{if } \tilde{n} \in \tilde{J}^\xi_1 \setminus \tilde{J}^\xi_2, \\
            0 & \text{if } \tilde{n} \in \tilde{J}^\xi_0 \setminus \tilde{J}^\xi_1.
        \end{cases}
    \end{equation*}
    \item\label{it: B5} For every $\alpha \leq \xi$, the pairs
    \begin{equation*}
        (\lambda^\alpha[J_2^\alpha], \lambda^\alpha[J_1^\alpha \setminus J_2^\alpha]) \text{ and } (\lambda^\alpha[J_1^\alpha], \lambda^\alpha[J_0^\alpha \setminus J_1^\alpha])
    \end{equation*}
    are not $(\mathfrak{B}(\xi+1 \setminus \{\alpha \}), \tilde{\mathfrak{B}}(\xi+1))$-separated; and the pairs
    \begin{equation*}
        (\gamma^\alpha[\tilde{J}_2^\alpha], \gamma^\alpha[\tilde{J}_1^\alpha \setminus \tilde{J}_2^\alpha]) \text{ and } (\gamma^\alpha[\tilde{J}_1^\alpha], \gamma^\alpha[\tilde{J}_0^\alpha \setminus \tilde{J}_1^\alpha])
    \end{equation*}
    are not $(\mathfrak{B}(\xi+1), \tilde{\mathfrak{B}}(\xi+1 \setminus \{\alpha\}))$-separated, where we recall that $\mathfrak{B}(\xi+1 \setminus \{\alpha \})$ denotes the Boolean subalgebra of $\mathcal{P}(\omega \times 2)$ given by
 \begin{equation*}
     \mathfrak{B}(\xi+1 \setminus \{\alpha \}) = [\{B_\eta^0, B_\eta^1: \eta \leq \xi, \eta \not = \alpha\} \cup \finset { \omega \times 2}],
 \end{equation*}
 and $\mathfrak{B}(\xi+1)$ denotes the Boolean subalgebra of $\mathcal{P}(\omega \times 2)$ given by
\begin{equation*}
    \mathfrak{B}(\xi+1) = [\{B_\eta^0, B_\eta^1: \eta \leq \xi\} \cup \finset { \omega \times 2}]. 
\end{equation*}
Analogous definitions apply to $\tilde{\mathfrak{B}}(\xi+1)$ and $\tilde{\mathfrak{B}}(\xi+1 \setminus \{ \alpha\})$. \\ 
 \end{enumerate}

We now describe how to carry out the construction for an ordinal $\xi <\mathfrak{c}$ assuming that we have already carried out the construction for all ordinals $\alpha < \xi$. That is, for every $\alpha < \xi$, we have chosen objects as specified in \ref{it: enu1}-\ref{it: enu5} so that conditions \ref{it: B1}-\ref{it: B5} hold. Note that this allow us to deal with both the initial and the induction steps at the same time. Furthermore, at this point of the induction process, the Boolean algebras $\mathfrak{B}(\xi)$ and $\tilde{\mathfrak{B}}(\xi)$ are defined. We show how to construct the objects at stage $\xi$.

Fix almost disjoint families $\mathcal{C}$ and $\tilde{\mathcal{C}}$, each of size $\mathfrak{c}$, consisting of infinite cylinders contained in $R_\xi \times 2$ and $\tilde{R}_\xi \times 2$, respectively. Since $(\mu^\xi_{n})_{n \in \omega}$, the $\ell_1(\omega \times 2)$-part of $(\lambda^\xi_{n})_{n \in \omega}$, is admissible, we can apply Lemma \ref{lmm: main} continuum many times—each time for a given cylinder $C = C_0 \times 2 \in \mathcal{C}$—to obtain sets $J_2(C) \subseteq J_1(C) \subseteq J_0(C) \subseteq C_0$ and $Z(C) \subseteq J_1(C) \times 2$, together with positive scalars $a(C)$ and $\delta(C)$ satisfying the conditions stated in the lemma. Furthermore, by Lemma \ref{lmm: only-countable-many} and Remark \ref{rem: uncountably-many-ways} we may assume that the pairs \begin{equation*}
        (\lambda^\xi[J_2(C)], \lambda^\xi[J_1(C) \setminus J_2(C)]) \text{ and } (\lambda^\xi[J_1(C)], \lambda^\xi[J_0(C) \setminus J_1(C)])
    \end{equation*}
    are not $(\mathfrak{B}(\xi), \tilde{\mathfrak{B}}(\xi))$-separated.

   Observe that for any $\alpha < \xi$ we have that the pairs
   \begin{equation*}
       (\lambda^\alpha[J_2^\alpha], \lambda^\alpha[J_1^\alpha \setminus J_2^\alpha]) \text{ and } (\lambda^\alpha[J_1^\alpha], \lambda^\alpha[J_0^\alpha \setminus J_1^\alpha])
   \end{equation*}
    are not $(\mathfrak{B}(\xi \setminus \{\alpha \}), \tilde{\mathfrak{B}}(\xi))$-separated; and the pairs
    \begin{equation*}
        (\gamma^\alpha[\tilde{J}_2^\alpha], \gamma^\alpha[\tilde{J}_1^\alpha \setminus \tilde{J}_2^\alpha]) \text{ and } (\gamma^\alpha[\tilde{J}_1^\alpha], \gamma^\alpha[\tilde{J}_0^\alpha \setminus \tilde{J}_1^\alpha])
    \end{equation*}
    are not $(\mathfrak{B}(\xi), \tilde{\mathfrak{B}}(\xi \setminus \{\alpha\}))$-separated. This follows from the fact that, by the induction hypothesis, \ref{it: B5} holds for any ordinal $\eta < \xi$. Indeed, if $\xi$ is a successor ordinal, say $\xi = \eta + 1$, then this is just \ref{it: B5} for $\eta$. Otherwise, if $\xi$ is a limit ordinal, this follows from the fact that an increasing family of non-separating algebras is also not separating.

   Therefore, we can apply Proposition \ref{prop: Haydon} to the almost disjoint family $\mathcal{Z} = \{Z(C): C \in \mathcal{C}\}$ and 
    \begin{itemize}[label=--]
        \item the family of algebra pairs $(\mathfrak{B}(\xi \setminus \{ \alpha\}), \tilde{\mathfrak{B}}(\xi))$ and the lists
        \begin{equation*}
             (\lambda^\alpha[J_2^\alpha], \lambda^\alpha[J_1^\alpha \setminus J_2^\alpha]) \text{ and } (\lambda^\alpha[J_1^\alpha], \lambda^\alpha[J_0^\alpha \setminus J_1^\alpha])
        \end{equation*}
        for $\alpha < \xi$.
    \item The family of algebra pairs $(\mathfrak{B}(\xi), \tilde{\mathfrak{B}}(\xi \setminus \{ \alpha\}))$ and the lists
    \begin{equation*}
        (\gamma^\alpha[\tilde{J}_2^\alpha], \gamma^\alpha[\tilde{J}_1^\alpha \setminus \tilde{J}_2^\alpha]) \text{ and } (\gamma^\alpha[\tilde{J}_1^\alpha], \gamma^\alpha[\tilde{J}_0^\alpha \setminus \tilde{J}_1^\alpha])
    \end{equation*}
    for $\alpha < \xi$.
    \end{itemize}  
    
    Thus, we can find $C \in \mathcal{C}$ such that:
    \begin{itemize}[label=--]
        \item for every $\alpha < \xi$, the pairs
            \begin{equation*}
                (\lambda^\alpha[J_2^\alpha], \lambda^\alpha[J_1^\alpha \setminus J_2^\alpha]) \hspace{5pt}\text{ and } \hspace{5pt} (\lambda^\alpha[J_1^\alpha], \lambda^\alpha[J_0^\alpha \setminus J_1^\alpha])
            \end{equation*}
         are not separated even by the larger algebra pair $(\mathfrak{B}(\xi \setminus \{ \alpha \})[Z(C)], \tilde{\mathfrak{B}}(\xi))$.
         \item For every $\alpha < \xi$, the pairs
                \begin{equation*}
                    (\gamma^\alpha[\tilde{J}_2^\alpha], \gamma^\alpha[\tilde{J}_1^\alpha \setminus \tilde{J}_2^\alpha]) \hspace{5pt}\text{ and } \hspace{5pt} (\gamma^\alpha[\tilde{J}_1^\alpha], \gamma^\alpha[\tilde{J}_0^\alpha \setminus \tilde{J}_1^\alpha])
                \end{equation*}
                 are not separated even by the larger algebra pair $(\mathfrak{B}(\xi)[Z(C)], \tilde{\mathfrak{B}}(\xi \setminus \{ \alpha \}))$. 
    \end{itemize}
     
     Having chosen $C \in \mathcal{C}$, we repeat the procedure in the second copy. Since $(\tilde{\mu}_{\tilde{n}})_{\tilde{n} \in \tilde{\omega}}$, the $\ell_1(\tilde{\omega} \times 2)$-part of $(\gamma^\xi_{\tilde{n}})_{\tilde{n} \in \tilde{\omega}}$, is $\sim$-admissible, we can apply Lemma \ref{lmm: main} continuum many times, each time for a cylinder $\tilde{C} = \tilde{C}_0 \times 2 \in \tilde{\mathcal{C}}$, to obtain sets $\tilde{J}_2(\tilde{C}) \subseteq \tilde{J}_1(\tilde{C}) \subseteq \tilde{J}_0(\tilde{C}) \subseteq \tilde{C}_0$ and $\tilde{Z}(\tilde{C}) \subseteq \tilde{J}_1(\tilde{C}) \times 2$, together with positive scalars $\tilde{a}(\tilde{C})$ and $\tilde{\delta}(\tilde{C})$ satisfying the conditions stated in the lemma. Again, by Lemma \ref{lmm: only-countable-many} and Remark \ref{rem: uncountably-many-ways} we may assume that the pairs 
    \begin{equation*}
        (\gamma^\xi[\tilde{J}_2(\tilde{C})], \gamma^\xi[\tilde{J}_1(\tilde{C}) \setminus \tilde{J}_2(\tilde{C})]) \text{ and } (\gamma^\xi[\tilde{J}_1(\tilde{C}))], \gamma^\xi[\tilde{J}_0(\tilde{C}) \setminus \tilde{J}_1(\tilde{C})])
    \end{equation*}
    are not $(\mathfrak{B}(\xi)[Z(C)], \tilde{\mathfrak{B}}(\xi))$-separated.
    
    We can invoke Proposition \ref{prop: Haydon}, this time for the almost disjoint family $\tilde{\mathcal{Z}} = \{\tilde{Z}[\tilde{C}]: \tilde{C} \in \tilde{\mathcal{C}} \}$ and
    \begin{itemize}[label = --]
        \item the family of algebra pairs $(\mathfrak{B}(\xi \setminus \{ \alpha \})[Z(C)], \tilde{\mathfrak{B}}(\xi))$ and the lists  
        \begin{equation*}
         (\lambda^\alpha[J_2^\alpha], \lambda^\alpha[J_1^\alpha \setminus J_2^\alpha]) \text{ and } (\lambda^\alpha[J_1^\alpha], \lambda^\alpha[J_0^\alpha \setminus J_1^\alpha])
    \end{equation*}
    for $\alpha < \xi$.
    \item The algebra pair $(\mathfrak{B}(\xi), \tilde{\mathfrak{B}}(\xi))$ and the pair
    \begin{equation*}
        (\lambda^\xi[J_2(C)], \lambda^\xi[J_1(C) \setminus J_2(C)]) \text{ and } (\lambda^\xi[J_1(C)], \lambda^\xi[J_0(C) \setminus J_1(C)]).
    \end{equation*}
    \item The family of algebra pairs $(\mathfrak{B}(\xi)[Z(C)], \tilde{\mathfrak{B}}(\xi))$ and the lists
    \begin{equation*}
        (\gamma^\alpha[\tilde{J}_2^\alpha], \gamma^\alpha[\tilde{J}_1^\alpha \setminus \tilde{J}_2^\alpha]) \hspace{5pt}\text{ and } \hspace{5pt} (\gamma^\alpha[\tilde{J}_1^\alpha], \gamma^\alpha[\tilde{J}_0^\alpha \setminus \tilde{J}_1^\alpha])
    \end{equation*}
    for $\alpha < \xi$.
    \end{itemize}

    Hence, we can find $\tilde{C} \in \tilde{\mathcal{C}}$ such that:
    \begin{itemize}[label=--]
        \item for every $\alpha < \xi$, the pairs
            \begin{equation*}
                (\lambda^\alpha[J_2^\alpha], \lambda^\alpha[J_1^\alpha \setminus J_2^\alpha]) \hspace{5pt}\text{ and } \hspace{5pt} (\lambda^\alpha[J_1^\alpha], \lambda^\alpha[J_0^\alpha \setminus J_1^\alpha])
            \end{equation*}
         are not separated even by the larger algebra pair $(\mathfrak{B}(\xi \setminus \{ \alpha \})[Z(C)], \tilde{\mathfrak{B}}(\xi)[\tilde{Z}(\tilde{C})])$.
         \item The pairs
             \begin{equation*}
                (\lambda^\xi[J_2(C)], \lambda^\xi[J_1(C) \setminus J_2(C)]) \text{ and } (\lambda^\xi[J_1(C)], \lambda^\xi[J_0(C) \setminus J_1(C)]).
            \end{equation*}
            are not separated even by the larger algebra pair $(\mathfrak{B}(\xi), \tilde{\mathfrak{B}}(\xi)[\tilde{Z}(\tilde{C})])$.
         \item For every $\alpha < \xi$, the pairs
                \begin{equation*}
                    (\gamma^\alpha[\tilde{J}_2^\alpha], \gamma^\alpha[\tilde{J}_1^\alpha \setminus \tilde{J}_2^\alpha]) \hspace{5pt}\text{ and } \hspace{5pt} (\gamma^\alpha[\tilde{J}_1^\alpha], \gamma^\alpha[\tilde{J}_0^\alpha \setminus \tilde{J}_1^\alpha])
                \end{equation*}
                 are not separated even by the larger algebra pair $(\mathfrak{B}(\xi)[Z(C)], \tilde{\mathfrak{B}}(\xi \setminus \{ \alpha \})[\tilde{Z}(\tilde{C})])$.
    \end{itemize}
    Finally, set
    \begin{enumerate}[label=(\roman*'), ref = (\roman*')]
        \item $J_2^\xi = J_2(C)$, $J_1^\xi = J_1(C)$, $J_0^\xi = J_0(C)$,
        \item $\tilde{J}_2^\xi = \tilde{J}_2(\tilde{C})$, $\tilde{J}_1^\xi = \tilde{J}_1(\tilde{C})$, $\tilde{J}_0^\xi = \tilde{J}_0(\tilde{C})$,
        \item $A_\xi = J_1(C) \times 2$, $B_\xi^0 = Z(C)$, $B_\xi^1 = A_\xi \setminus B_\xi^0$; $\tilde{A}_\xi = \tilde{J}_1(\tilde{C}) \times 2$, $\tilde{B}_\xi^0 = \tilde{Z}(\tilde{C})$, $\tilde{B}_\xi^1 = \tilde{A}_\xi \setminus \tilde{B}_\xi^0$,
        \item $a_\xi = a(C)$, $\tilde{a}_\xi = \tilde{a}(\tilde{C})$, $\delta_\xi = \delta(C)$ and $\tilde{\delta}_\xi = \tilde{\delta}(\tilde{C})$.
    \end{enumerate}
    By construction, conditions \ref{it: B1}--\ref{it: B5} are satisfied. 
    Indeed, \ref{it: B1}--\ref{it: B2} are precisely \eqref{eq: A1} in Lemma~\ref{lmm: main}, which hold by our choices in the construction. 
    Moreover, \ref{it: B5} follows directly from the construction, since we ensured an appropriate choice guaranteeing separation of measures. 
    This completes the induction. \\

With the construction completed according to the inductive procedure, we now define the almost disjoint families
    \begin{equation*}
    \mathcal{A} = \{A_\xi : \xi < \mathfrak{c}\}, \hspace{3pt}
    \tilde{\mathcal{A}} = \{\tilde{A}_\xi : \xi < \mathfrak{c}\}, 
    \hspace{3pt} \mathcal{B} = \{B_\xi^0, B_\xi^1 : \xi < \mathfrak{c}\},
    \text{ and } \tilde{\mathcal{B}} = \{\tilde{B}_\xi^0, \tilde{B}_\xi^1 : \xi < \mathfrak{c}\}.
    \end{equation*}
These families will ensure that the corresponding $X$ and $\tilde{X}$ are not isomorphic to Banach lattices, as we will demonstrate in the next section.

\bigskip
\section{Proof of Theorem \ref{th: main-lattices}} \label{sec: lattice}

Let $\mathcal{A}, \tilde{\mathcal{A}}, \mathcal{B}$ and $\tilde{\mathcal{B}}$ be the almost disjoint families built in Section \ref{sec: proof-main-theorem}. We show that, for these families, the spaces $X$ and $\tilde{X}$ such that
\begin{equation*}
    \JL(\mathcal{B}, \tilde{\mathcal{B}}) = X \oplus \tilde{X,}
\end{equation*}
defined as in Subsection \ref{subsec: two-spaces}, are not isomorphic to Banach lattices. Since the proofs for $X$ and $\tilde{X}$ are identical, we only do it for the former. 

Observe that $X$ is a $1$-complemented subspace of $\JL(\mathcal{B}, \tilde{\mathcal{B}})$. Therefore $X^*$ is $1$-complemented in $\JL(\mathcal{B}, \tilde{\mathcal{B}})^* \equiv \ell_1(K_\mathcal{B} \cup K_{\tilde{\mathcal{B}}})$ and hence it is linearly isometric to $\ell_1(\Gamma)$ for some $\Gamma$. Furthermore, observe that
\begin{equation*}
    \{\delta_{(n,0}), \delta_{(n,1)}: n \in \omega\} \cup \{\delta_{(\tilde{n},0}), \delta_{(\tilde{n},1)}: \tilde{n} \in \tilde{\omega}\}
\end{equation*}
is a $1$-norming set for $\JL(\mathcal{B}, \tilde{\mathcal{B}})$ and thus, in particular, $X$ has a countable $1$-norming set. That is, the conditions of Proposition \ref{lmm: non-banach-lattice-criteria} are satisfied.

For convenience, we also state the following consequences of the construction. Recall that $\{(\lambda_n^\xi)_{n \in \omega}: \xi < \mathfrak{c}\}$ is an ordering of $\Phi$, as described in the construction, and that we have a decomposition $\lambda_n^\xi = (\mu_n^\xi + \overline{\nu}_n^\xi) + \tilde{\nu}_n^\xi$.

\begin{enumerate}[label = (P\arabic*), ref = (P\arabic*),leftmargin=*]
        \item \label{it: p1} For every $\xi < \mathfrak{c}$ let $b_\xi = \inf_{n \in \omega} |\mu_n^\xi(\{(n, 0) \})| > 0$. By \ref{it: B1}, there exist  $a_\xi \geq b_\xi$ and $0 < \delta_\xi < b_\xi/11$ such that
        \begin{equation*}
        \mu^\xi_n(B_\xi^0) - \mu^\xi_n(B_\xi^1) \approx_{3\delta_\xi} \begin{cases}
            2a_\xi & \text{if } n \in J^\xi_2,  \\
            -2a_\xi & \text{if } n \in J^\xi_1 \setminus J^\xi_2, \\
            0 & \text{if } n \in J^\xi_0 \setminus J^\xi_1.
        \end{cases}
        \end{equation*}
        \item \label{it: p3} For any $\xi < \mathfrak{c}$, the pairs
        \begin{equation*}
            (\lambda^\xi[J_2^\xi], \lambda^\xi[J_1^\xi \setminus J_2^\xi]) \text{ and } (\lambda^\xi[J_1^\xi], \lambda^\xi[J_0^\xi \setminus J_1^\xi])
            \end{equation*}
        are not $(\mathfrak{B}(\mathfrak{c} \setminus \{\xi\}), \tilde{\mathfrak{B}})$-separated. This follows since \ref{it: B5} holds for every $\eta < \mathfrak{c}$, and the union of an increasing family of non-separating algebras is again not separating.
\end{enumerate}

We are now ready for the proof.

\begin{proof}[Proof of Theorem \ref{th: main-lattices}]
    By Proposition \ref{lmm: non-banach-lattice-criteria} and the previous discussion, it is enough to show that $X$ has the \eqref{eq: DP}. Thus, we start by fixing a norming sequence $(e^*_n)_{n \in \omega}$ for $X$. Our objective is to find $f \in X$ such that no $g \in X$ satisfies
    \begin{equation*}
        \langle e^*_n, g \rangle = |\langle e_n^*, f \rangle| \hspace{5pt} \text{ for every } n \in \omega.
    \end{equation*}

    Since $(e_n^*)_{n \in \omega}$ is norming for $X$ and $f_k = \mathds{1}_{(k,0)} - \mathds{1}_{(k,1)} \in Y \subseteq X$ (see Subsection \ref{subsec: two-spaces} for the definition of $Y$), passing to a subsequence of $(e_n^*)_{n \in \omega}$, we may assume that
    \begin{equation*}
        0 < \inf_{n \in \omega} |\langle e_n^*, f_n \rangle|.
    \end{equation*}
    
    From Remarks \ref{rem: norming-as-subsets-measures} and \ref{rem: decomp-of-measures}, we can write $e_n^* = (\mu_n + \overline{\nu}_n) + \tilde{\nu}_n$, so that $\mu_n(\{(k, 0)\}) + \mu_n(\{(k, 1)\}) = 0$ for every $k \in \omega$, $\overline{\nu}_n(B_\alpha^0) + \overline{\nu}_n(B_\alpha^1) = 0$ for every $\alpha < \mathfrak{c}$ and $\norm{e^*_n}_1 = \norm{\mu_n}_1 + \norm{\overline{\nu}_n}_1 + \norm{\tilde{\nu}_n}_1 \leq 1$. In other words, the coding of $e_n^*$ is of \ref{type-1}.
    
    Observe that
    \begin{equation*}
        0 <\frac{1}{2}\inf_{n \in \omega} |\langle e_n^*, f_n \rangle| = \inf_{n \in \omega} |e_n^*(\{(n,0)\})| = \inf_{n \in \omega} |\mu_n(\{(n,0)\})| =: b.
    \end{equation*}
    From the above discussion and the construction in Section \ref{sec: proof-main-theorem}, we can find $\xi < \mathfrak{c}$ such that $(e_n^*)_{n \in \omega} = (\lambda^\xi_n)_{n \in \omega}$. For convenience, we drop the $\xi$ from the measures, and simply write $(e_n^*)_{n \in \omega} = (\lambda_n)_{n \in \omega} = (\mu_n + \overline{\nu}_n + \tilde{\nu}_n)_{n \in \omega}$. Similarly, we also denote $\delta \equiv \delta_\xi$ and $a \equiv  a_\xi $ according to \ref{it: p1}.
    
    Recall that, from the way the enumeration was carried out, for any $n \in \omega$ and $i = 0, 1$ we have $\overline{\nu}_n(B_\alpha^i) = 0$ whenever $\alpha \geq \xi$. We claim that
    \begin{equation*}
        f = \mathds{1}_{B^0_\xi} - \mathds{1}_{B^1_\xi} \in Y \subseteq X
    \end{equation*}
    witnesses \eqref{eq: DP}, which will finish the proof. Assume the contrary, so that we can find $g \in X$ such that
    \begin{equation*}
        \langle\lambda_n, g\rangle = |\langle\lambda_n, f\rangle |  \hspace{5pt} \text{ for every } n \in \omega.
    \end{equation*}
    Note that $\langle\overline{\nu}_n, f\rangle = 0$ since $\overline{\nu}_n(B_\xi^i) = 0$ for every $n \in \omega$ and $i = 0,1$. Clearly, $\langle\tilde{\nu}_n, f \rangle = 0$ and thus $\langle\lambda_n, f\rangle = \langle\mu_n, f\rangle$.

    Since the subspace
    \begin{equation*}
        \{h + \tilde{h}: h \in Y \text{ simple } \mathfrak{B}\text{-measurable, and } \tilde{h} \in \JL(\tilde{A}) \text{ simple } \tilde{\mathfrak{B}}\text{-measurable} \}
    \end{equation*}
    is dense in $X$, it follows that we can find a simple $\mathfrak{B}$-measurable $h \in Y$ and a simple $\tilde{\mathfrak{B}}$-measurable function $\tilde{h} \in \JL(\tilde{\mathcal{A}})$ such that $\norm{g - (h + \tilde{h})} < \delta$. Hence
    \begin{equation*}
         |\langle\mu_n, f \rangle| = |\langle\lambda_n, f\rangle | = \langle \lambda_n, g\rangle \approx_\delta \langle\lambda_n, h + \tilde{h}\rangle, \hspace{5pt} \text{ for every } n \in \omega.
    \end{equation*}

    Without loss of generality, we may assume that $h = rf + s$ where $r \in \mathbb{R}$ and $s \in Y$ is a simple $\mathfrak{B}(\mathfrak{c \setminus \{ \xi\}})$-measurable function, so that
    \begin{equation}\label{eq: 1}
        |\langle\mu_n, f \rangle| \approx_\delta r \langle\mu_n, f\rangle + \langle\lambda_n, s+ \tilde{h}\rangle, \hspace{5pt} \text{ for every } n \in \omega.
    \end{equation}
    Furthermore, we may assume that $r \geq 0$, the case $r < 0$ follows in the same way. Observe that \ref{it: p1} gives
    \begin{equation}\label{eq: 2}
        \langle\mu_n, f\rangle = \mu_n(B_\xi^0) - \mu_n(B_\xi^1) \approx_{3\delta} \begin{cases}
            2a & \text{if } n \in J_2^\xi,  \\
            -2a & \text{if } n \in J^\xi_1 \setminus J_2^\xi, \\
            0 & \text{if } n \in J_0^\xi \setminus J^\xi_1,
        \end{cases}
    \end{equation} 
    and thus in particular
    \begin{equation}\label{eq: 3}
        |\langle\mu_n, f\rangle| \approx_{3\delta}  \begin{cases}
            2a & \text{if } n \in J_1^\xi,  \\
            0 & \text{if } n \in J_0^\xi \setminus J^\xi_1.
        \end{cases}
    \end{equation}
    Combining \eqref{eq: 1}, \eqref{eq: 2} and \eqref{eq: 3} yields
    \begin{equation}\label{eq: 4}
        \begin{cases}
            2a  \approx_{\delta(4r + 3)} 2ra + \langle\lambda_n, s + \tilde{h}\rangle &  \text{if } n \in J_2^\xi, \\
            2a  \approx_{\delta(4r + 3)} -2ra + \langle\lambda_n, s + \tilde{h}\rangle& \text{if } n \in J^\xi_1 \setminus J_2^\xi, \\
            0 \approx_{\delta(4r + 3)}   \langle\lambda_n, s + \tilde{h}\rangle& \text{if } n \in J_0^\xi \setminus J^\xi_1.
        \end{cases}
    \end{equation}

    Suppose first that $0 \leq r \leq 1/2$, the first two equations in \eqref{eq: 4} give
    \begin{equation*}
        \langle\lambda_n, s + \tilde{h} \rangle \geq 2(1 - r)a - \delta(4 + 3r) \geq a - \frac{11}{2}\delta, \hspace{5pt} \text{ for every } n \in J_1^\xi,
    \end{equation*}
   while the third equation in \eqref{eq: 4} yields
    \begin{equation*}
        \langle\lambda_n, s + \tilde{h}\rangle \leq \delta(4r + 3) \leq \frac{11}{2} \delta , \hspace{5pt} \text{ for every } n \in J_0^\xi \setminus J^\xi_1.
    \end{equation*}
    Therefore, since $a \geq b >  11\delta$, we get
    \begin{equation*}
        \langle\lambda_n, s + \tilde{h} \rangle - \langle\lambda_k, s + \tilde{h}\rangle \geq a - 11\delta > 0
    \end{equation*}
    for any $n \in J_1^\xi$ and $k \in J_0^\xi \setminus J_1^\xi$. Lemma \ref{lmm: simple-functions-separate-sets} now gives that the pair
    \begin{equation*}
        \{\lambda_n: n \in J_1^\xi\} = \lambda_n[J_1^\xi] \hspace{5pt} \text{ and } \hspace{5pt} \{\lambda_n: n \in J_0^\xi \setminus J_1^\xi\} = \lambda_n[J_0^\xi \setminus J_1^\xi]
    \end{equation*}
    is $(\mathfrak{B}(\mathfrak{c}\setminus \{\xi \}), \tilde{\mathfrak{B}})$-separated, which contradicts \ref{it: p3}.

    Suppose now that $r > 1/2$, then by \eqref{eq: 4} we deduce that
    \begin{align*}
    \langle \lambda_n, s + \tilde{h} \rangle - \langle \lambda_k, s + \tilde{h} \rangle 
    &\geq 2a(1+r) - \delta(4+3r) - \bigl(2a(1-r) + \delta(4+3r)\bigr) \\
    &= 2r(2a - 3\delta) - 8\delta \geq 2a - 11\delta > 0.
    \end{align*}
    for any $n \in J_1^\xi \setminus J_2^\xi$ and $k \in  J_2^\xi$. By Lemma \ref{lmm: simple-functions-separate-sets} we have that the pair
    \begin{equation*}
        \{\lambda_n: n \in J_2^\xi\} = \lambda_n[J_2^\xi]  \hspace{5pt} \text{ and } \hspace{5pt} \{\lambda_n: n \in J_1^\xi \setminus J^\xi_2\} = \lambda_n[J_1^\xi \setminus J^\xi_2]
    \end{equation*}
    is $(\mathfrak{B}(\mathfrak{c}\setminus \{\xi \}), \tilde{\mathfrak{B}})$-separated, which again contradicts \ref{it: p3}. This finishes the proof.
\end{proof}

\noindent\textbf{Acknowledgements.} This paper forms part of the author’s PhD research at Lancaster University, conducted under the supervision of Professor N. J. Laustsen. The author is deeply grateful to Professor Laustsen for helpful comments and suggestions in the presentation of the manuscript. He is also thankful to Matthew Daws and Pedro Tradacete for their comments on a preliminary version of this manuscript.

He acknowledges with thanks the funding from the EPSRC (grant number \break EP/W524438/1) that has supported his studies.

% \bibliographystyle{ieeetr.bst}
% \bibliography{refs}

\end{document}